\documentclass[conference]{IEEEtran}

\usepackage{amsmath,epsfig}
\usepackage{cite,subfigure}
\usepackage[normalem]{ulem}
\usepackage[T1]{fontenc}
\usepackage[utf8]{inputenc}
\usepackage[bookmarks=false]{hyperref}
\usepackage{amssymb}
\usepackage{upgreek}
\usepackage{amsthm}
\usepackage{enumerate}
\usepackage{color}

\newtheorem{theorem}{Theorem}[section]
\newtheorem*{theorem*}{Theorem}

\newtheorem{define}[theorem]{Definition}

\DeclareMathOperator{\diag}{\bf diag}

\newcommand{\beq}{\begin{equation}}
\newcommand{\eeq}{\end{equation}}
\newcommand{\beqn}{\begin{equation*}}
\newcommand{\eeqn}{\end{equation*}}
\newcommand{\beqa}{\begin{equation}\begin{aligned}}
\newcommand{\eeqa}{\end{aligned}\end{equation}}
\newcommand{\beqna}{\begin{equation*}\begin{aligned}}
\newcommand{\eeqna}{\end{aligned}\end{equation*}}
\newcommand{\ii}{\item}
\newcommand{\enum}{\begin{enumerate}}
\newcommand{\enuma}{\begin{enumerate}[(a)]}
\newcommand{\eenum}{\end{enumerate}}
\newcommand{\norm}[1]{||#1||}
\newcommand{\inprod}[1]{\langle #1 \rangle}

\newcommand{\lp}{\left(}
\newcommand{\rp}{\right)}

\newcommand{\vphi}{\varphi}

\newcommand{\reals}{\mathbb{R}}

\newcommand{\integers}{\mathbb{Z}}

\newcommand{\dual}[1]{\widetilde{#1}}

\newcommand{\abs}[1]{\left|#1\right|}

\newcommand{\bolx}{\boldsymbol{x}}

\newcommand{\bolf}{\boldsymbol{f}}

\newcommand{\bolb}{\boldsymbol{\beta}}

\newcommand{\bol}[1]{\boldsymbol{#1}}

\newcommand{\twonorm}[1]{\norm{#1}_2}

\title{Using Frame Theoretic Convolutional Gridding for Robust Synthetic Aperture Sonar Imaging}
\author{\IEEEauthorblockN{John McKay$^\blacklozenge$, Anne Gelb$^\boxminus$, Vishal Monga$^\blacklozenge$, Raghu G. Raj$^\odot$}
\IEEEauthorblockA{
$^\blacklozenge$Department of Electrical Engineering, 
Pennsylvania State University, 
State College, PA\thanks{McKay, Monga, \& Raj were supported by ONR Grant 0401531.}\\
$^\boxminus$Department of Mathematics, 
Dartmouth College, 
Hanover, NH\thanks{A. Gelb was partially supported by AFOSR FA9550-15-1-0152}\\
$^\odot$U.S. Naval Research Laboratory, Washington, DC}
}
\IEEEoverridecommandlockouts 
\begin{document}
\maketitle

\begin{abstract}
Recent progress in synthetic aperture sonar (SAS) technology and processing has led to significant advances in underwater imaging, outperforming previously common approaches in both accuracy and efficiency. There are, however, inherent limitations to current SAS reconstruction methodology.  In particular, popular and efficient Fourier domain SAS methods require a $2D$ interpolation which is often ill conditioned and inaccurate, inevitably reducing robustness with regard to speckle and inaccurate sound-speed estimation.  To overcome these issues, we propose using the {\em frame theoretic convolution gridding} (FTCG) algorithm to handle the non-uniform Fourier data.  FTCG extends upon non-uniform fast Fourier transform (NUFFT) algorithms {by casting the NUFFT as an approximation problem given Fourier frame data.}  The FTCG has been show to yield improved accuracy at little more computational cost.  Using simulated data, we outline how the FTCG can be used to enhance current SAS processing.
\end{abstract}

\section{Introduction}
\label{sec:intro}

Synthetic aperture Sonar has proven to be one of the most important developments underwater imaging.  Indeed, {SAS} can produce high resolution images of much larger areas than was previously possible \cite{hansen2011introduction}.  Most {reconstruction} algorithms developed for Sonar fall into two categories: back-projection and  Fourier domain (FD).  Back-projection methods have shown promise in terms of resolution but requires significant computational effort.  On the other hand, FD methods ``can reconstruct images in real time on standard workstation computers'' \cite{hayes2009synthetic} but with varying quality  \cite{jakowatz2006comparison}.

Many different methods fall into the FD category, including wavenumber and range Doppler algorithms in SAS \cite{gough2005fast} and polar-format algorithms (PFAs) in SAR \cite{fan2015improved} \cite{zhu2007range} \cite{fan2014polar}.  In each, a $2D$ interpolation step emerges as {the main issue limiting} resolution.  From a numerical analysis point of view, this problem is even more difficult since no one optimal mathematical framework exists for two dimensional interpolation \cite{hayes2009synthetic}.  Of the ways this interpolation step has been addressed, Strolt mapping \cite{gough1997unified}, chirp z-transforms \cite{zhu2007range}, and non-uniform fast Fourier Transforms (NUFFTs) \cite{andersson2012fast,fan2015improved} have had some success in SAS.  Although chirp z-transforms seemlingly yield the fastest and most accurate algorithms, they are susceptible to irregular vehicle movements and therefore suitable in only the most ideal cases \cite{zhu2007range}.  On the other hand, NUFFTs - namely fast Gaussian gridding (FGG) NUFFTs - have demonstrated a robustness to noise and jitter in the Fourier space in synthetic aperture radar applications while still retaining reasonable computational efficiency  \cite{greengard2004accelerating,andersson2012fast,fan2014polar}.  It is with this in mind that we present a similar but further refined algorithm for tomographic imaging:  frame theoretic convolutional gridding (FTCG).  The FTCG algorithm was developed in \cite{gelb2014frame} where the authors demonstrate how by cleverly incorporating frame theoretic concepts into NUFFTs, high quality image reconstructions can be had at comparable computational efficiencies and with properties guaranteeing bounded error.  

In this paper, we will look to briefly introduce SAS imaging processing, acquaint the reader to the FTCG algorithm in the one dimensional manner presented in \cite{gelb2014frame}, show how it can be adapted to two dimensions as to be applied to SAS, and then detail applications of this method to simulated data.

\section{Synthetic Aperture Sonar}
\label{sec:SASintro}

\begin{figure}[!t]\centering
\includegraphics[width=1\columnwidth]{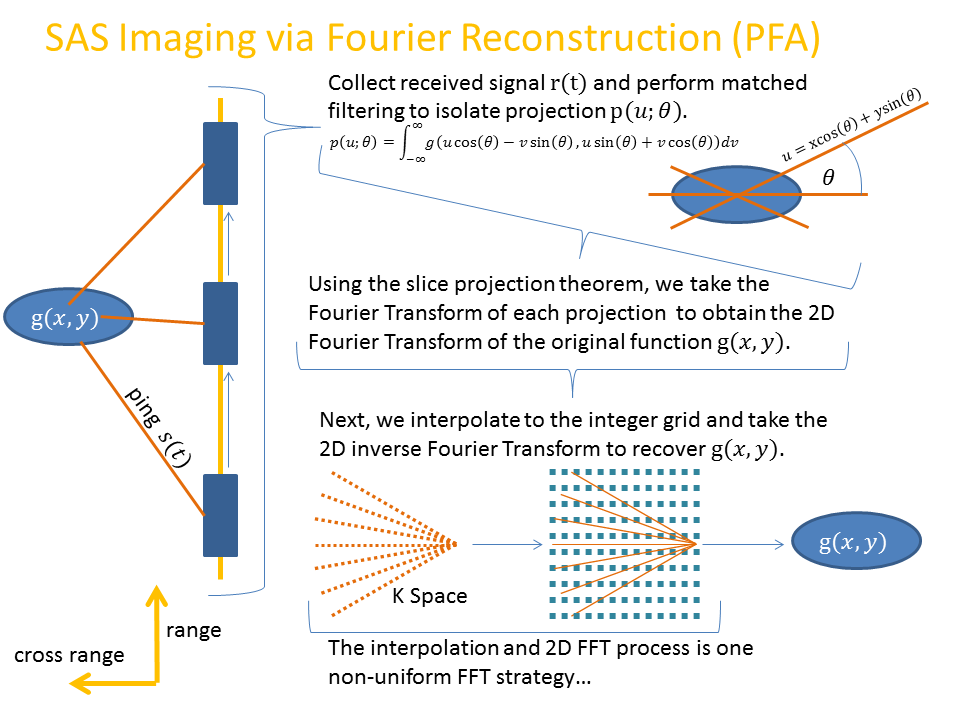}
\caption{SAS imaging procedure with focus on Fourier domain methods for reconstruction.}\label{SAS}
\end{figure}

SAS is as powerful as it is because it is able to overcome unavoidable hardware constraints.   The along track resolution of Sonar, Radar, and other tomographic systems is proportional to the width of the aperture of the signal-emitting device \cite{hayes2009synthetic}.  Unfortunately, though, there is a trade-off between narrow apertures and device size: the smaller the aperture, the larger the necessary array of sensors.  This factor has limited Sonar imaging in the past, but as illustrated in Figure \ref{SAS}, SAS processing has been able to work by utilizing a series of successive pings and combining them as to mimic a smaller aperture.  {Of the two methods we mentioned previously, note that in this setting back-projection approaches essentially perform several costly inverse Radon transforms and FD methods leverage the Fourier slice projection theorem (FSPT) to attempt to reconstruct scenes from frequency data quickly \cite{soumekh1999synthetic}.  Herein lies the problem:  the FSPT will \emph{rarely} (if ever) yield data points in some sort of equispaced arrangement.  This means that one must employ some sort of correction or use a flexible strategy like a $2D$ NUFFT \cite{soumekh1999synthetic,breed2001k}.  If a correction (i.e. interpolation) is used, it may struggle with data collected under irregular vehicle movements and/or incorrect sound-speed estimation which can incur unsightly artifacts like blurring \cite{cook2009analysis}}.

{This all means that a well designed, robust NUFFT algorithm can be quite valuable for SAS imaging, let alone other tomographic fields.}  Some success has been had in synthetic aperture radar in applying NUFFTs to their similar FD algorithms and tend to offer a far greater robustness than interpolation-based approaches \cite{soumekh1999synthetic}.  \cite{gough2005fast} provided an introductory view of NUFFTs in SAS, but this is not a well-studied problem.

\section{Image Reconstruction from Non-Uniform Fourier Data}
{Below we present a  description of standard NUFFT and the FTCG algorithm that is used to enhance its performance.}

\subsection{The NUFFT}

The NUFFT essentially maps irregular Fourier data to a uniform grid by attempting to quantify the ``contribution'' of Fourier basis element to each data point.  For example, if we consider an image as a (one-dimensional) reflectivity function $f$ from sound waves, then its NUFFT approximation $A_{cg}(f)$ is given by
\beqa\label{eq:ffind}
A_{cg}(f)= \sum\limits_{\abs{n} \leq N}\sum\limits_{{\abs{m-\lambda_n}\leq K}}\alpha_n\hat{f}(\lambda_n)\hat{w}(m-\lambda_n)\tfrac{e^{2\pi i m x}}{w(x)},
\eeqa
where $f$ is assumed to be piece-wise smooth on $[0,1]$, $\hat f(\lambda_n), n=-N,\dots,N$, are its Fourier coefficients on some raster $\{\lambda_n\}$, {and $K$ is a truncation threshold determined by the support of $\hat{w}$.}  {The calculation performed in (\ref{eq:ffind}) can be viewed as  $g_K(x)/w(x)$ where $g_K$ is the $K$th order Fourier partial sum approximation of $g = fw$.  Here  $w$ is a positive window function that is {\em essentially} compactly supported  in the domain $[0,1]$ while also having the property that its Fourier coefficients, $\hat{w}$ are essentially compactly supported.  This is necessary to reduce the impact of aliasing while still maintaining computational efficiency.  The support in the Fourier domain of $\hat{w}$ determines the threshold $K$ in (\ref{eq:ffind}). Finally, $\hat{w}$ should have an analytical expression so as not to require calculation on point values $(m-\lambda_n)$.}  Gaussian window functions have shown promise in SAR \cite{greengard2004accelerating}, though there are several other options.  

Much work has been done in determining the quadrature weights $\alpha_m$ employed in approximating  the convolutional integral $\hat{g} = \hat f\ast \hat w$.  In particular, iterative methods have been demonstrated to be effective \cite{pipe1999sampling}.

For the purpose of our discussion, we will rewrite $A_{cg}$ in condensed form. {It is convenient to replace the sum limits $|m - \lambda_n| \le K$ with $|m|\le M$ which does not affect the solution (only zeroes are added in) but makes comparison to the new method more straightforward.  We write}
\beqa\label{eq:cgmatrix}
A_{cg}(f)(x) = \sum\limits_{\abs{m}\leq M}\gamma_m\psi_m{(x)}\text{ where } \bol{\gamma}=\Omega D \hat\bolf,
\eeqa
where
\begin{eqnarray}
\Omega&=&[\hat w(m-\lambda_n)]_{{\abs{m}\leq M},\abs{n}\leq N},\nonumber\\
D&=&\diag(\alpha_{-N}\,\,\cdots\,\,\alpha_N),\nonumber\\
\hat\bolf &=&\begin{pmatrix} \hat f(\lambda_{-N})\,\,\cdots\,\,\hat f(\lambda_N)\end{pmatrix}^T,
\label{eq:D_alpha}
\end{eqnarray}
and $ \psi_m(x)=\exp(2\pi imx)/w(x)$, which will be especially relevant later when we discuss frames.

As we have alluded to before, NUFFTs are known to have manageable computational cost.  While the exact cost depends on the window function, the idea is that $\hat w$ has small support in the Fourier domain.  In this way, the convolutional expense is small (which is represented by the matrix $\Omega$) so that $\bol{\gamma}$ can be calculated quickly and the subsequent reconstruction can be done by the weighted sums represented in \eqref{eq:cgmatrix}.

\subsection{Frame Approximation}

{The {\em frame theoretic} FFT was designed in \cite{song2013approximating} to improve the accuracy of $A_{cg}$ by incorporating the properties of frames into the design of the quadrature weights $\alpha$.  A brief review of the method is provided below.}

\begin{define}
A frame is a set $\{\vphi_j\}$ in a Hilbert space $\mathcal{H}$ such that
\beqa\label{eq:frame}
A\norm{h}^2_\mathcal{H} \leq \sum\limits_{n\in {\mathbb{Z}}}\abs{ \inprod{h,\vphi_n} }\leq B\norm{h}_\mathcal{H}\,\,\,\,\,\forall h\in\mathcal{H}
\eeqa
for some $A,\,B>0$.  Any such set of vectors that satisfy \eqref{eq:frame} \emph{must} span $\mathcal{H}$ even if the set contains redundancies.  
\end{define}
Observe that a frame is similar to a basis for $\mathcal{H}$ but with a relaxation on linear independence. 

Although the set $\{\vphi_n=\exp(2\pi i\lambda_n)\}$, $m \in \mathbb{Z}$  is typically not a basis if $\lambda_n \ne n$,  for some sampling patterns $\vphi_n$ forms a frame \cite{benedetto2002mri}. It was also demonstrated in  \cite{gelb2014frame} that even if this is not the case, the assumption that $\{\vphi_n\}$ is a frame can still be quite useful.  

The {\em frame operator} is defined as
\beqa\label{eq:S}
Sf = \sum\limits_{\abs{n}\leq \mathbb{Z}}\inprod{f,\vphi_n}\vphi_n,
\eeqa
where $\hat{f}_n = \inprod{f,\vphi_n}$ are known as the {\em frame coefficients}.  Observe that, since $S$ is both linear and invertible \cite{christensen2003introduction}, we can say
\beqa\label{eq:framereco}
f=S^{-1}Sf = S^{-1}\sum\limits_{\abs{n}\leq \mathbb{Z}}\hat{f}_n\vphi_n
= \sum\limits_{\abs{n}\leq \mathbb{Z}}\hat{f}_nS^{-1}\vphi_n.
\eeqa
{Typically the dual frame $\dual\vphi_n= S^{-1}\vphi_n$ is not known in closed form, so an approximation is needed.  Moreover, in applications we only know a finite number of frame coefficients, e.g.~  $\hat{f}_n = \hat{f}(\lambda_n)$, $n = -N,\cdots,N$.  Hence we are interested in the partial sum frame approximation}
\beqa\label{eq:framereco_finite}
A_{fr}(f) = \sum\limits_{\abs{n}\leq N}\hat{f}_n\dual\vphi_n.
\eeqa
The  method developed  in \cite{song2013approximating} approximates the dual frame $\dual\vphi_n$ so that (\ref{eq:framereco_finite}) converges to (\ref{eq:framereco}).  The method is based on projecting $\vphi$ onto a so called {\em admissible frame} for $\mathcal{H}$.  The convergenge of (\ref{eq:framereco_finite}) is then dependent on the properties of the admissible frame. 
\begin{define}
\label{def:admissible}
An admissible frame $\{\psi_m\}_{\abs{m}\leq M}$ of $\{\phi_i\}_{\abs{i}\leq N}$ is a frame chosen so that
\enum[i]
\ii  $\displaystyle \abs{\inprod{\psi_i,\psi_k}}\leq c(1+\abs{i-k})^{-t}$ for some $c>0$, $t>1$, and all $i,k$.
\ii $\displaystyle \abs{\inprod{\phi_i,\psi_k}}\leq d(1+\abs{i-k})^{-s}$ for some $d>0$, $s>1/2$, and all $i,k$.
\eenum
\end{define}
Given an admissible frame, an approximation to the  dual frame can then be written as
\beqa\label{eq:sumof}
\dual\vphi_n = \sum\limits_{\abs{m}\leq M}b_{n,m}\psi_m.
\eeqa
In \cite{song2013approximating} it was shown that the coefficients in (\ref{eq:sumof}) can be solved for using  $B=[b_{n,m}]_{\abs{n}\leq N,\abs{m}\leq M}$ with
\beqa\label{eq:binv}
B=\Psi^\dagger\text{ for }\Psi=[\inprod{\phi_n,\psi_m}]_{\abs{n}\leq N,\abs{m}\leq M},
\eeqa
where $\Psi^\dagger$ is the Moore Penrose pseudo-inverse.  From the discussion above, we can now approximate any $h\in\mathcal{H}$ as
\beqa\label{eq:tie}
h \approx \sum\limits_{\abs{m}\leq M}\sum\limits_{\abs{n}\leq N}\inprod{h,\vphi_n}b_{n,m}\psi_m
\eeqa

To link these ideas to (\ref{eq:cgmatrix}), we define 
$\psi_m=\exp(2\pi i m x)/w(x)$, which was also used in \cite{gelb2014frame} with respect to this data frame.   Thus, we modify the finite frame approximation in (\ref{eq:framereco_finite}) to our problem of recovering reflectivity function $f$ as\footnote{We note that (\ref{eq:framereco_finite}) is not equivalent to (\ref{eq:afr}), but we use the same variable so as to avoid cumbersome notation.}
\beqa\label{eq:afr}
A_{fr}(f)=\sum\limits_{\abs{m}\leq M}\sum\limits_{\abs{n}\leq N} b_{n,m}\hat{f}(\lambda_n)\tfrac{\exp(2\pi i  m x)}{w(x)}
\eeqa
or, in matrix form
\beqa\label{eq:afrm}
A_{fr}(f)=\sum\limits_{\abs{m}\leq M}\beta_m\psi_m \text{ where } \bolb=B\hat{\bolf}.
\eeqa

The frame approximation has several advantages over NUFFTs, not the least of which is better reconstruction as well as a known upper error bound for certain rasters, even when corrupted with noise \cite{gelb2014frame}.  Thus, the frame approximation provides a certain degree of confidence that other methods cannot.  On the other hand, this comes with considerable cost; $A_{fr}$ requires {$\mathcal{O}(MN)$ operations, in addition to the (off-line) calculation of $B$,} making it infeasible for many machines.

\subsection{Frame Theoretic Convolutional Gridding (FTCG)}

The FTCG strategy bridges the gap between $A_{cg}$ and $A_{fr}$ that takes advantage of the efficiency of the NUFFT but provides  higher accuracy by replacing the traditional vector of quadrature weights with a {\em matrix} of quadrature weights informed by frame theory.  Specifically, the FTCG approximation $A(f)$ is given by
\beqa\label{eq:ftcg}
A(f) = \sum\limits_{\abs{m}\leq M} {\tau}_m\psi_m \text{ for }\bol{\tau}= \Omega C \hat\bolf,
\eeqa
where $C$ is an $2r-1$-banded matrix.  The comparison to the traditional convolutional gridding method is easily seen in the matrix formula for $A_{cg}$ in (\ref{eq:cgmatrix}) and in particular $D$ in (\ref{eq:D_alpha}).  Here $C$, derived below,  is banded as compared to diagonal $D$.\footnote{{As discussed in \cite{gelb2014frame}, for $D$ diagonal, the coefficients of $A_{cg}$ are determined using standard numerical quadrature.  In this case, the convergence of $A_{cg}$ inherently depends upon the density of data points.  Specifically, the density of points must be inversely proportional to the number of data points.  However, in SAS imaging, this is not necessarily the case, as increasing $N$ can lengthen the spectrum.}  
{Hence in order to guarantee convergence, a different convolution approximation is needed.}}

We first describe how to design $D$, the diagonal matrix of quadrature weights $\{\alpha_n\}_{n = -N}^N$.   Observe that the Frobenious norm difference between $A_{cg}$ and $A_{fr}$ is
\beqa\label{eq:dif}
\norm{A_{cg}-A_{fr}}_F\leq K\twonorm{\bol\gamma - \bolb},
\eeqa
where $K$ is a constant stemming from frame bounds, and $\beta$ and $\gamma$ are defined in (\ref{eq:afrm}) and (\ref{eq:cgmatrix}) respectively.  {Furthermore, we see that}
\beqa\label{eq:dif2}
\twonorm{\bol\gamma-\bolb}\leq\norm{\lp\Psi^H\Psi\rp^{-1}}_F\norm{\Psi^H\Psi\Omega D - \Psi^H}\twonorm{\hat\bolf},
\eeqa
which can be simplified to show
\beqa\label{eq:together}
\twonorm{A_{cg}-A_{fr}}\leq \norm{(\Psi^H\Psi)^{-1}}_F\twonorm{\hat\bolf}\norm{\Psi\Omega D -I}_F.
\eeqa 
Therefore, solving
\beqa\label{eq:difopt}
\min\limits_{D} \norm{\Psi\Omega D - I}_F,
\eeqa
provides a straight-forward way to find optimal diagonal elements to solve \eqref{eq:difopt}.  Let us now remove the constraint for $D$ to be diagonal, then we see that the minimum is achieved for
\beqa\label{eq:dif3}
D=\lp\Psi\Omega\rp^\dagger.
\eeqa
Inserting this into \eqref{eq:cgmatrix} yields coefficients
\beqa\label{eq:insert}
{\bol\gamma}= \Omega D \hat\bolf = \Omega\Omega^\dagger\Psi^\dagger\bolf = \Psi^\dagger\bolf=\bol\beta,
\eeqa
{so that by (\ref{eq:afrm})} we end up with $A_{fr}$!  
Therefore, we have a link between the two methods and we can formally designate a trade-off between computational efficiency and accuracy based on the band we impose for the inverted matrix.  In other words, the coefficient matrix $C$ in FTCG is such that
\beqa\label{eq:cr}
C = (\Psi\Omega\odot B_{r})^\dagger
\eeqa
where $\odot$ is the Hadamard product and $B_{r}$ is a binary matrix with ones along a $r-1\times r-1$ band and zeros elsewhere.  This isolation of the elements of $\Psi\Omega$ along a band provides the key aspect to the FTCG method:  increasing $r$ improves the accuracy of the method while reducing its efficiency.  In \cite{gelb2014frame} it was suggested to choose $r\approx \log(2N+1)$ so that the cost of the convolutional sums mirrors that of subsequent FFTs, making for a minimal increase in computation over $A_{cg}$ but with the promise of a better reconstruction.

\begin{figure}[!b]\centering
\includegraphics[width=.9\columnwidth]{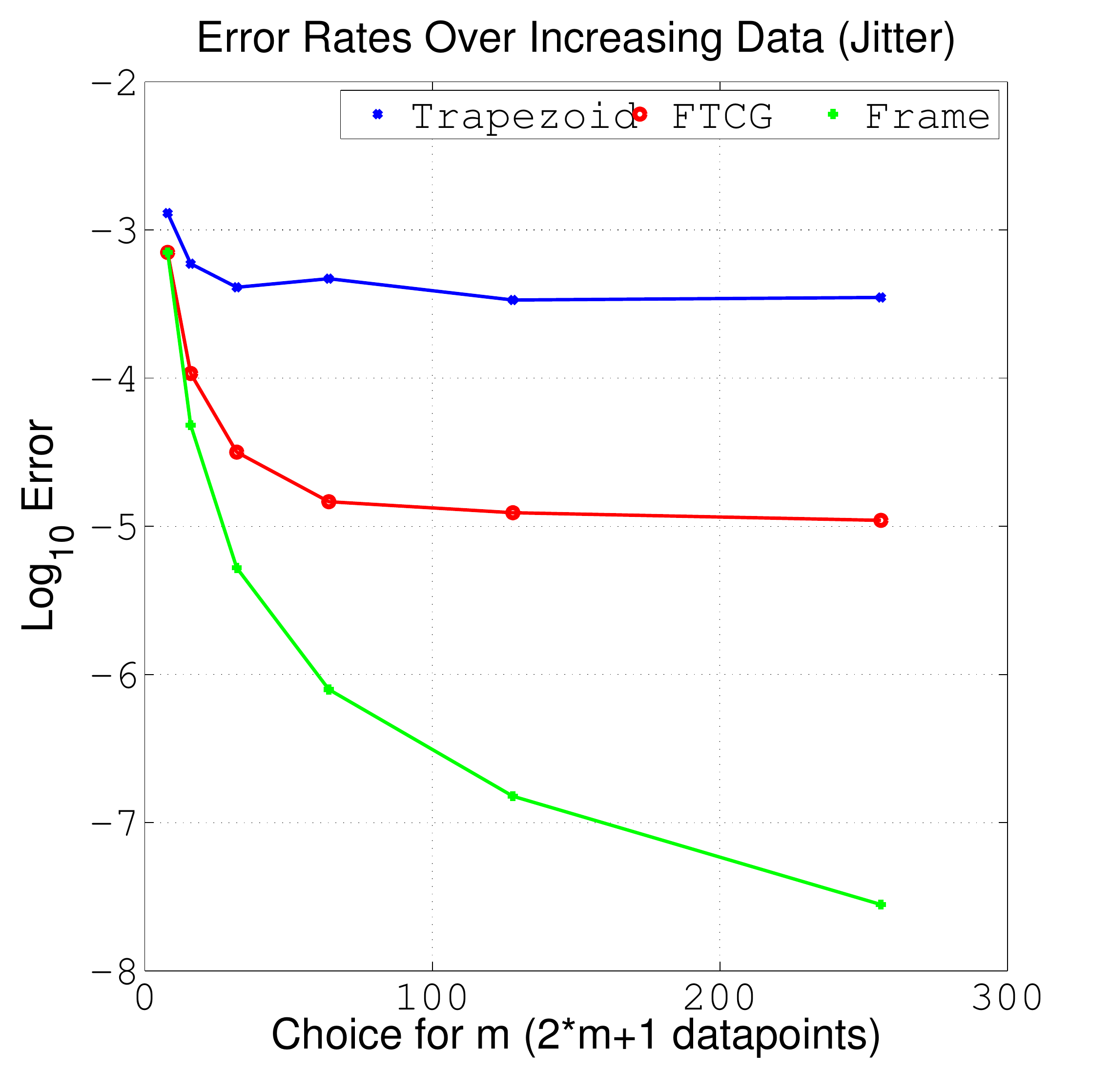}
\caption{$\ell_2$ error of 3 1-D methods over varied raster sizes.}\label{errorComp}
\end{figure}

Another important advantage of the FTCG over the standard CG method is improved convegence properties.  As Figure \ref{errorComp} shows in a $1D$ experiment over a jittered raster, simple convolutional gridding garners little benefit from more data. In fact, as was observed in \cite{viswanathan2010reconstruction}, increasing the number of Fourier samples potentially {\em reduces} the accuracy  since {it amounts to increasing the size of the spectral domain.  Hence the density between points actually grows. As discussed previously, the convergence of numerical quadrature methods is linked to the density of data points, i.e. $\max{|\omega_{n+1}-\omega_n|}$, for $n = -N,\cdots,N-1$.}  On the other hand, the convergence of the frame approximation {is based on the  {\em truncation error}}, that is, the accuracy improves with increased $N$.  However, it requires  computing the inverse of an $N \times N$ full matrix.  The FTCG provides a balance between the two.   The accuracy increases with bandwidth $r$, allowing practitioners with large datasets to adjust the band for both efficiency and accuracy.

\section{Two Dimensional FTCG}

The authors of \cite{song2015two} present a $2D$ frame approximation for non-uniform Fourier data points.  The method is based on an extension of the one-dimensional admissible frame methodology into $d > 1$ dimensions.  We briefly describe the method below for $d = 2$.

Suppose we are given a $2D$ raster of data where each element,  $\bol\lambda_{\bol m}\in\reals^2$, is such that for some $M_1, M_2\in\mathcal{N}$
\beqa
\bol m\in\mathcal{I}_m=\{[m_1,m_2]\in\integers^2\,:\,\abs{m_1}\leq M_1,\,\abs{m_2}\leq M_2\}.
\eeqa
In sequel we will refer to $\bol n\in\mathcal{I}_n$ defined similarly for some $N_1,N_2\in\mathcal{N}$.

The data frames $\{\vphi_{\bol m}\}_{\bol m\in\mathcal{I}_m}$ are defined as
\beqa\label{eq:twoddataframe}
\vphi_{\bol m}(\bolx)=\exp(2\pi j \inprod{\bol\lambda_{\bol m},\bolx}).
\eeqa
As in the $1D$ case, we define an {\em admissible frame} for $\vphi_{\bol m}$. 
\begin{define} $\{\psi_{\bol n}\}$, $\bol n\in\mathcal{I}_n$, is an admissible frame with respect to $\{\vphi_{\bol m}\}$ in two dimensions given that the following are satisfied
\enum[a)]
\ii\label{admone}  $\abs{\inprod{\psi_{\bol n},\psi_{\bol \ell}}}\leq \gamma_0(1+\twonorm{\bol n - \bol\ell})^{-t}$\\ for $\gamma_0>0$, $t>2$, $\bol n,\bol\ell\in\mathcal{I}_n$
\ii\label{admtwo}  $\abs{\inprod{\vphi_{\bol m},\psi_{\bol n}}}\leq \gamma_1(1+\abs{m_1-n_1})^{-s}(1+\abs{m_2-n_2})^{-s}$\\ for $\gamma_1,s>0$, $\bol m\in\mathcal{I}_m$, and $\bol n\in\mathcal{I}_n$
\eenum
\end{define}
For example,
\beqa\label{eq:twodadmissible}
\psi_{\bol n}(\bolx) = \exp(2\pi j\inprod{\bol n,\bolx}),\hspace{.2in}\bol n\in\mathcal{I}_n,
\eeqa
satisfies the admissibility requirements.

Although the admissible frame construction yields accurate reconstruction for a variety of different data arrangements with bounded maximal error, it is computationally inefficient.  Further, $\psi_{\bol n}$ as suggested in \eqref{eq:twodadmissible} is not readily applicable to a $2D$ NUFFT scheme.  Hence we seek to adapt the method in \cite{song2015two} towards an accurate and efficient  $2D$ FTCG algorithm.

\begin{figure*}\centering
\includegraphics[width=.32\textwidth]{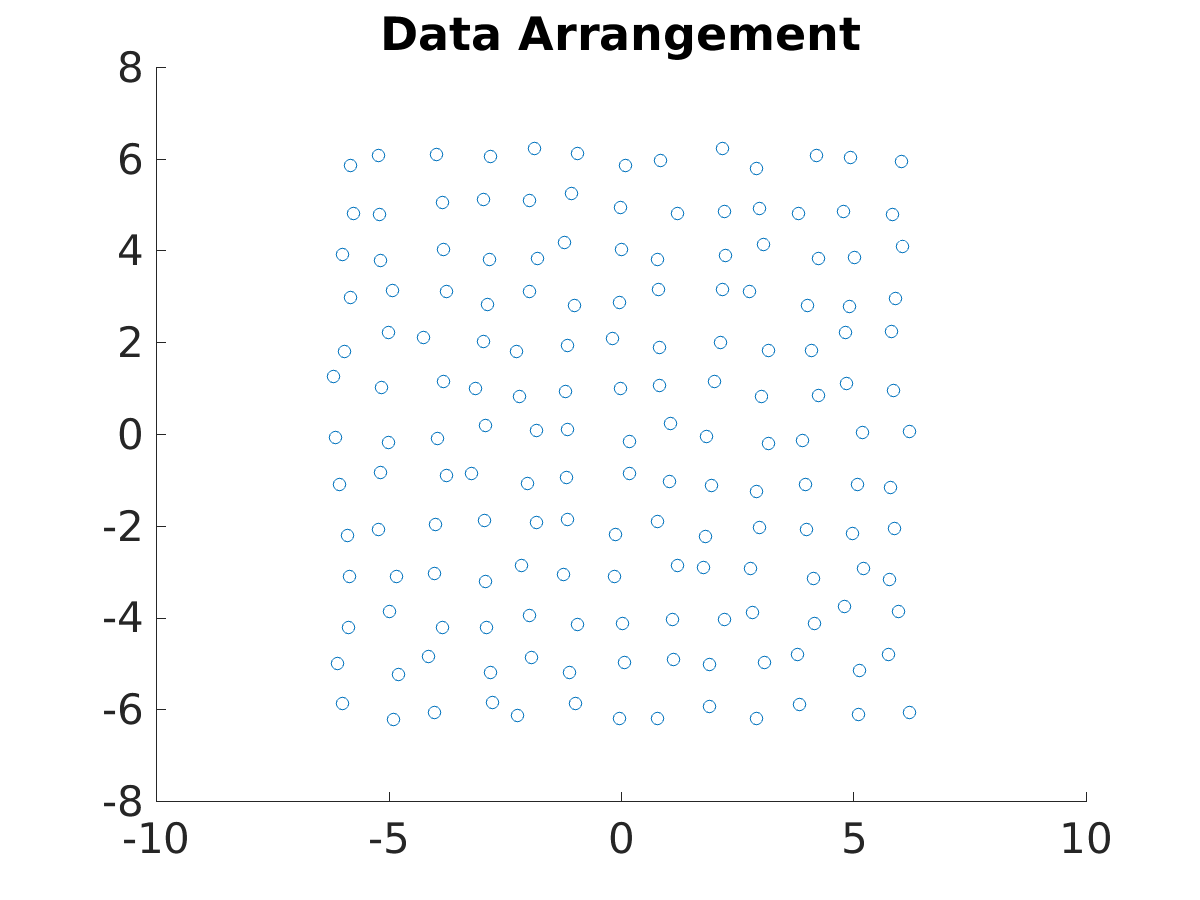}
\includegraphics[width=.32\textwidth]{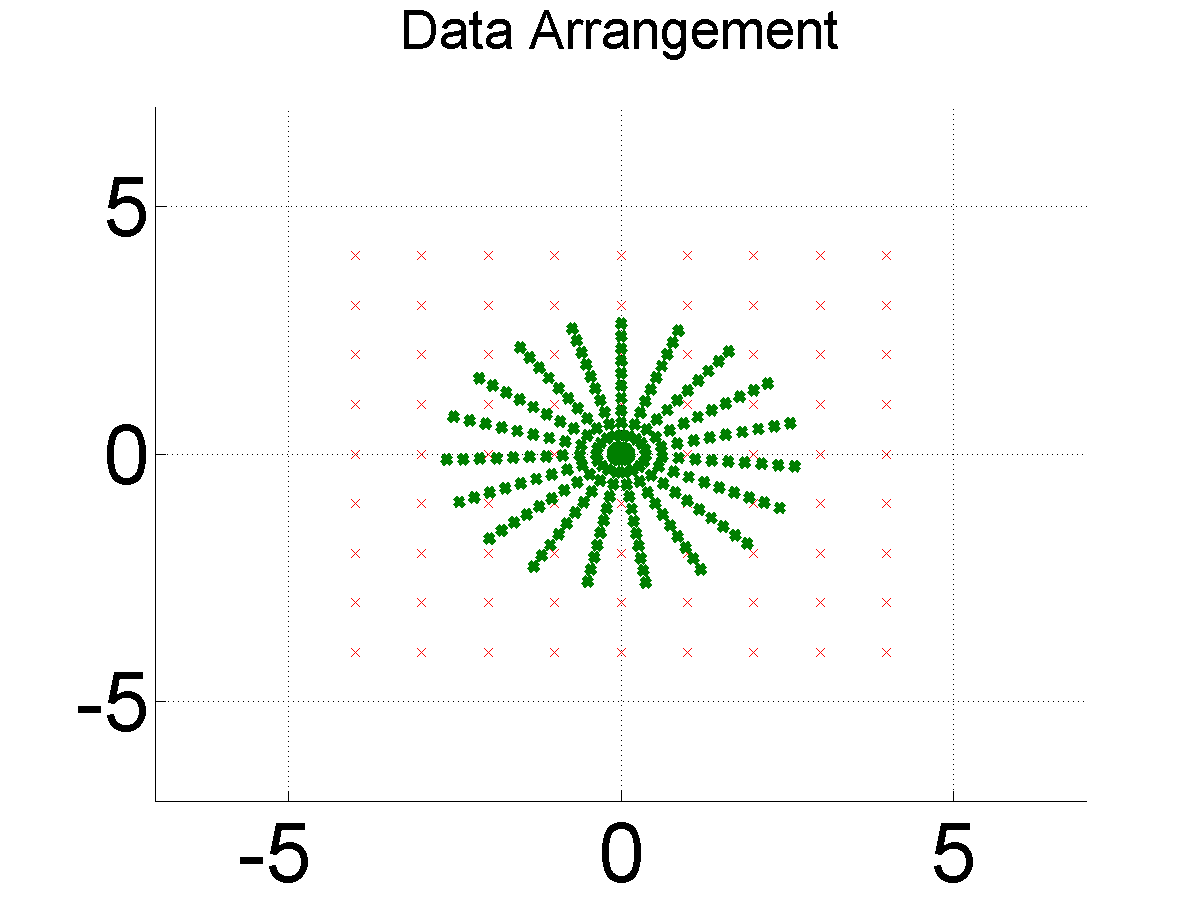}
\includegraphics[width=.32\textwidth]{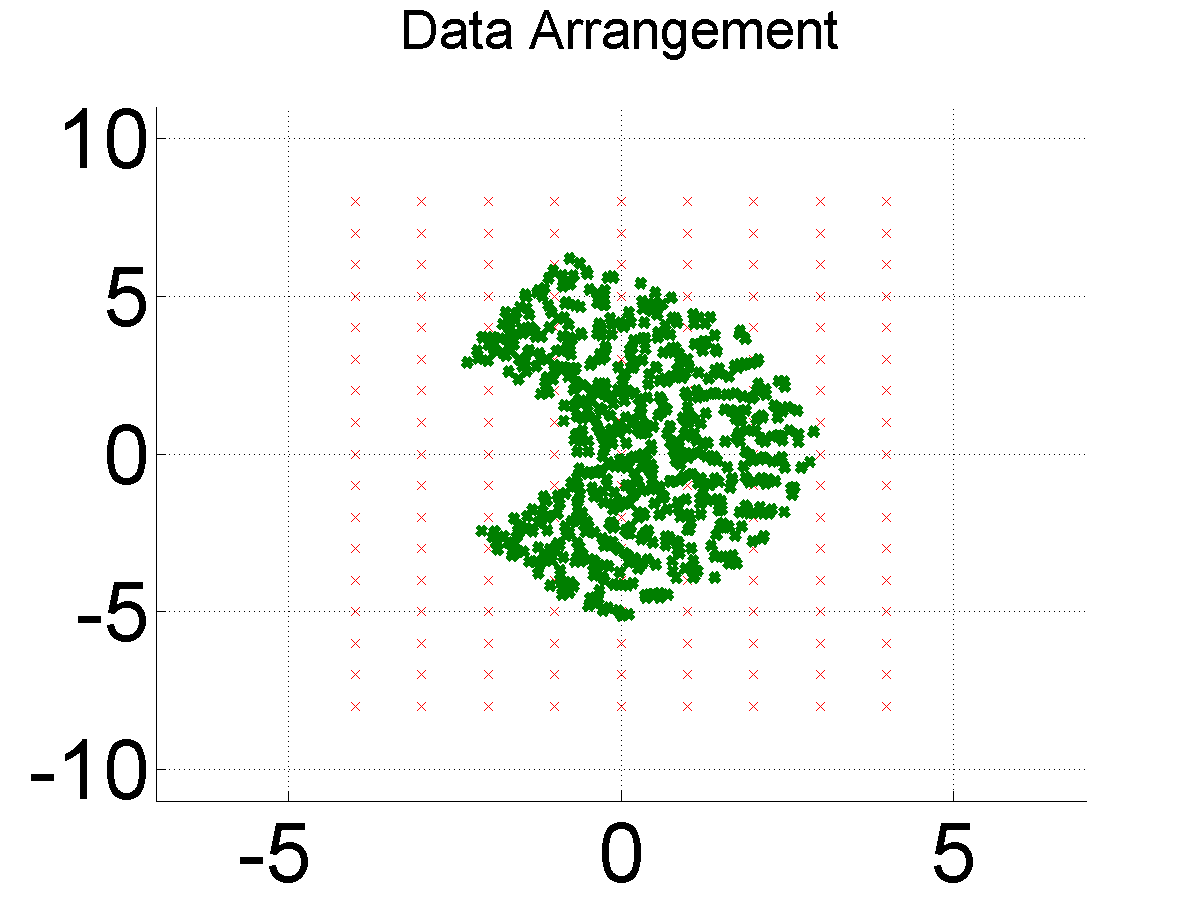}
\caption{Example Fourier data raster arrangements of noisy grids (left), asterisk patterns (center), and SAS-like set-ups (right).}\label{fig:jitgrid}
\end{figure*}

As in the $1D$ case, we use the convolutional gridding method for construction.  In this regard, let us define $w(\bolx)$ as a bounded $2D$ window function with essentially compact support over $[0,1]^2$ (that is, $w$ decays rapidly outside of that region) and $\{\zeta_{\bol n}\}_{\bol n\in\mathcal{I}_n}$ be such that
\beqa\label{eq:newadm}
\zeta_{\bol n}(\bolx) =\tfrac{1}{w(\bolx)}\psi_{\bol n}(\bolx).
\eeqa
For $\psi_{\bol n}$ as defined in \eqref{eq:twodadmissible}, we can show that $\{\zeta_{\bol n}\}_{\bol n\in\mathcal{I}_n}$ is also an admissible frame for $\{\phi_{\bol m}\}_{\bol m\in\mathcal{I}_m}$.  To do so, consider the two criteria for $2D$ admissibility:
\beqna
\text{a)  }  &\abs{\inprod{\zeta_{\bol n},\zeta_{\bol \ell}}}&=&\abs{\inprod{\tfrac{\psi_{\bol n}}{w(\bolx)},\tfrac{\psi_{\bol \ell}}{w(\bolx)}}}\\
&&=&\tfrac{1}{w(\bolx)^2}\abs{\inprod{\psi_{\bol n},\psi_{\bol \ell}}}\\
&&\leq & \tfrac{1}{B}\abs{\inprod{\psi_{\bol n},\psi_{\bol \ell}}}\\
&&\leq & \tfrac{\gamma_0}{B}(1+\twonorm{\bol n - \bol\ell})^{-t},
\eeqna
using sufficient $\gamma_0$, $t$ from the admissibility of $\{\psi_{\bol n}\}$ and bound $B>0$ for $w$, and
\beqna
\text{b)  }\abs{\inprod{\vphi_{\bol m},\zeta_{\bol n}}}&=\tfrac{1}{B}\abs{\inprod{\vphi_{\bol m},\psi_{\bol n}}}\\
&\leq\tfrac{\gamma_1}{B}(1+\abs{m_1-n_1})^{-s}(1+\abs{m_2-n_2})^{-s},
\eeqna
again utilizing $\gamma_1$ and $s$ from the admissibility of $\psi_{\bol n}$.  Thus, having satisfied the criteria \ref{admone}) and \ref{admtwo}), we see that, indeed, $\{\zeta_{\bol n}\}_{\bol n\in\mathcal{I}_n}$ is an admissible frame for our data frame.

The $2D$ FTCG is analogous to the $1D$ version.  Assuming that our original function $f:\reals^2\to\reals^2$ can be reconstructed from a separable Fourier space, we augment \eqref{eq:cgmatrix} to\footnote{{For comparative purposes we write the sum to include all values $\bol\ell\in\mathcal{I}_L$. As before, when computing the sum, only a subset of those values are used since $\hat{w}$ is compactly supported.}}
\beqa\label{eq:twodcg}
A^2_{cg}(f)=\sum\limits_{\bol\ell\in\mathcal{I}_L}\sum\limits_{\bol m\in\mathcal{I}_m} \alpha_{\bol m}\hat f(\bol\lambda_{\bol m})\hat w(\bol\ell -\bol\lambda_{\bol m})\tfrac{\exp(2\pi j \inprod{\bol\ell,\bolx})}{w(\bolx)},
\eeqa
with truncation term 
\beqna
\bol\ell\in\mathcal{I}_L=\{[\ell_1,\ell_2]\in\integers^2:\,\abs{\ell_1}\leq L_1,\,\abs{\ell_2}\leq L_2\}.
\eeqna
A more compact form for \eqref{eq:twodcg} is expressed as
\beqa\label{eq:twodcgmatrix}
A^2_{cg}(f)&=\sum\limits_{\bol m\in\mathcal{I}_m}\gamma_{\bol m}\psi_{\bol m},
\eeqa
where $\gamma_{\bol m}=\Omega_2 D_2\hat\bolf_2$ and each component defined as
\begin{align}
\Omega_2&=[\hat w(\bol\ell - \bol\lambda_{\bol m})]_{\bol m\in\mathcal{I}_m,\bol\ell\in\mathcal{I}_L}\\
D_2&=\diag(\{\alpha_{\bol m}\})\\
\hat\bolf_2 &= [\hat f(\lambda)_{\bol m})]^T.
\end{align}
The $2D$ frame approximation with $\{\zeta_{\bol n}\}$ is given by
\beqa\label{eq:twodframe}
A^2_{fr}(f)=\sum\limits_{\bol n\in\mathcal{I}_n} \beta_{\bol n}\zeta_{\bol n},
\eeqa
such that for $\bol\beta=[\beta_{\bol n}]_{\bol n\in\mathcal{I}_n}^T$ we have
\beqa\label{eq:twodframeB}
\bol\beta = \Psi_2^\dagger\hat\bolf_2 \text{ where }\Psi_2=[\inprod{\vphi_{\bol m},\zeta_{\bol n}}]_{\bol m\in\mathcal{I}_m,\bol n\in\mathcal{I}_n}.
\eeqa

As in the $1D$ case, the relationship between $A^2_{cg}$ and $A^2_{fr}$ can be understood by the optimization problem expressed by \eqref{eq:difopt} with the $2D$ analogous terms.  In particular, we can construct the $2D$ FTCG approximation of band $r$ as
\beqa\label{eq:twodftcg}
A^2(f) = \sum\limits_{\bol m\in\mathcal{I}_m}\tau_{\bol m}\zeta_{\bol m} \text{ for } \bol{\tau}=\Omega_2 C_2\hat\bolf_2,
\eeqa
where $C_2=(\Psi_2\Omega_2\odot B_r)^\dagger$.

We note here that this construction assumes separability of the Fourier space.  Future investigations will consider the non-separable case.

\section{Experiments}

To show the effectiveness of the FTCG algorithm, we present the following results using synthetic data mimicking Sonar image capture.  {Future investigations  will include real SAS data.}

We look to show the viability of FTCG in Sonar by demonstrating its ability on a noisy grid, an asterisk pattern, and then a more challenging side-scan-SAS-like arrangement.  Each of these arrangements were chosen to reflect an inherent robustness to data sampling patterns. {Indeed, as demonstrated below, the FTCG is readily used for a variety of data environments.}

For each case, we looked to reconstruct the function
\beqa\label{eq:eq2}
f(\bolx)=\sin(4\pi x_1)\sin(2\pi x_2).
\eeqa
which is similar to that used in \cite{song2015two}.  

\subsection{Noisy Grid Data Arrangement}

\begin{figure*}\centering
\begin{minipage}{1\textwidth}
\includegraphics[width=.5\textwidth]{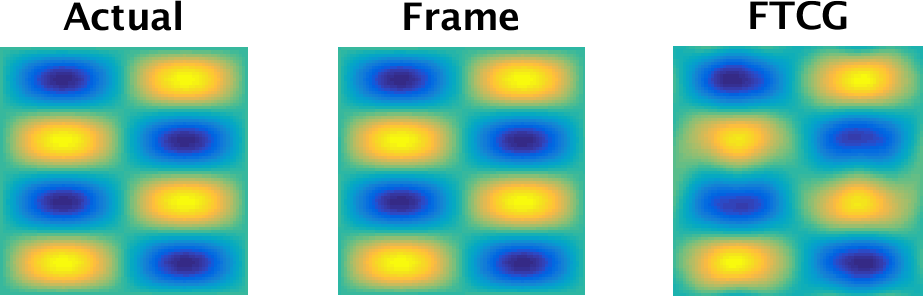}
\includegraphics[width=.245\textwidth]{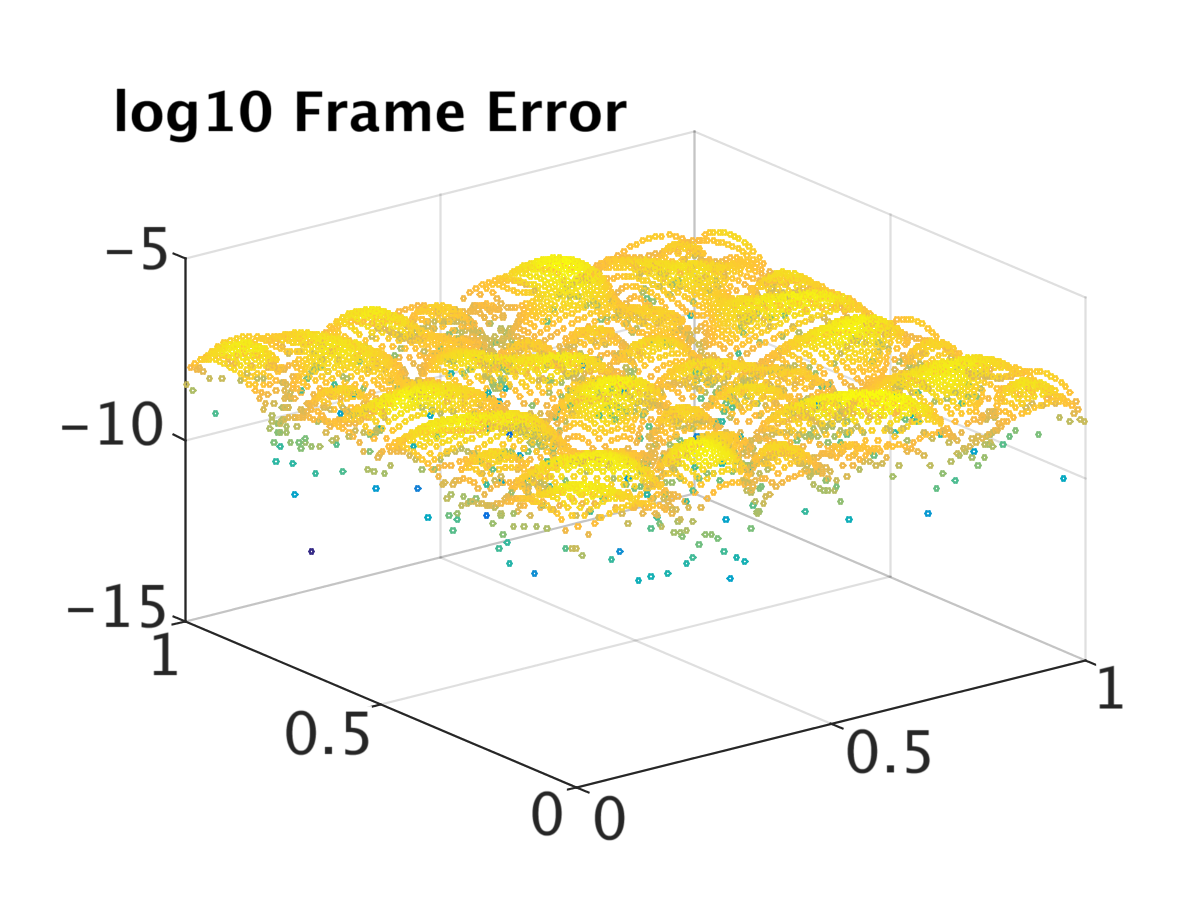}
\includegraphics[width=.245\textwidth]{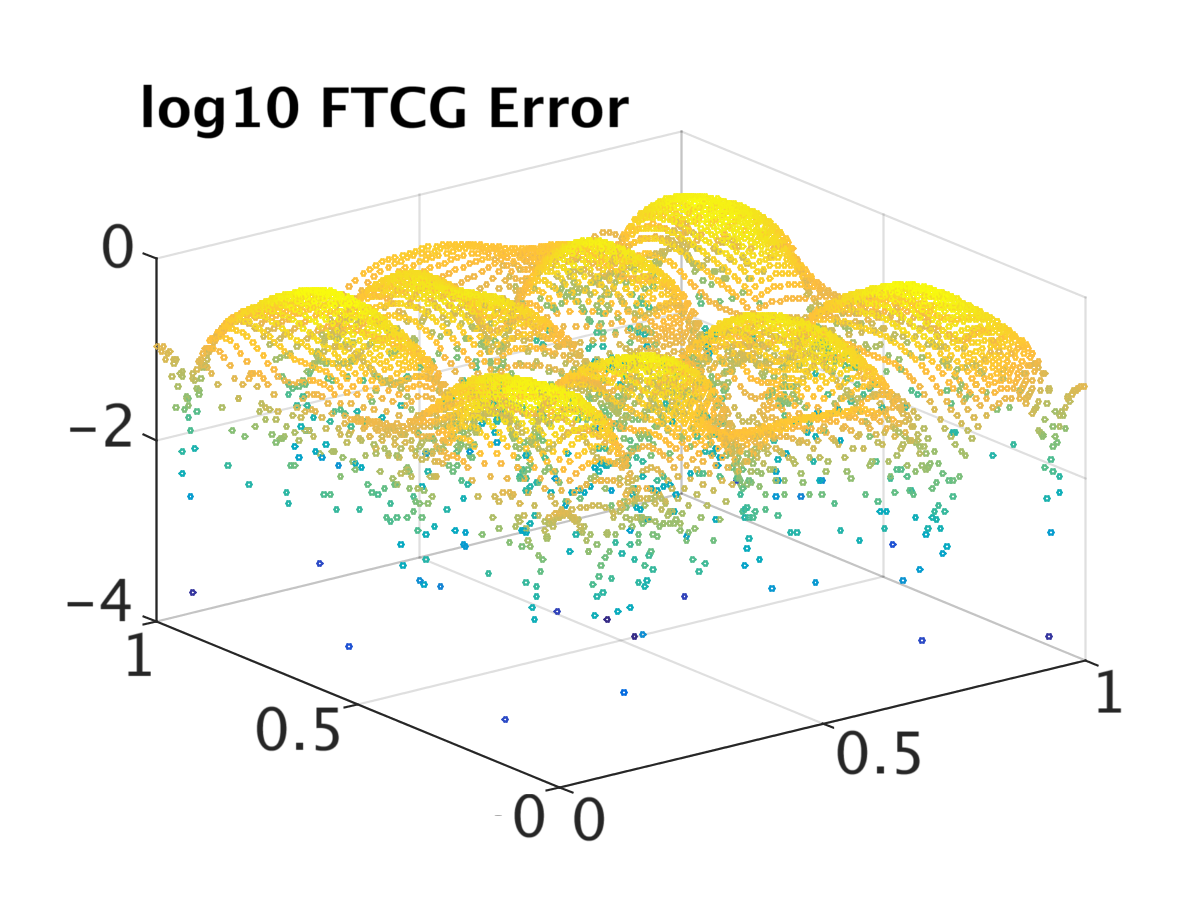}
\caption{Function reconstruction from noisy data arrangement.}\label{fig:nois}
\end{minipage}

\begin{minipage}{1\textwidth}
\includegraphics[width=.5\textwidth]{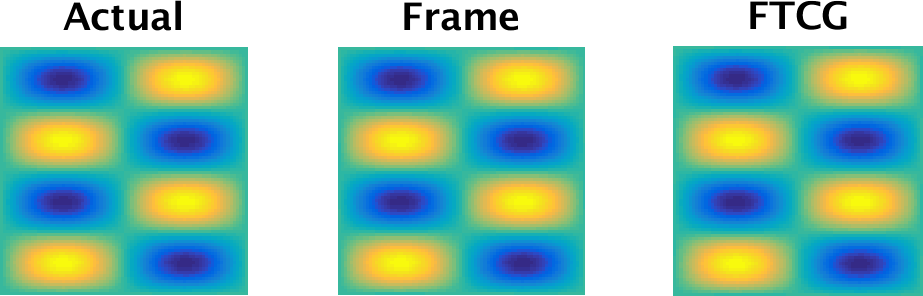}
\includegraphics[width=.245\textwidth]{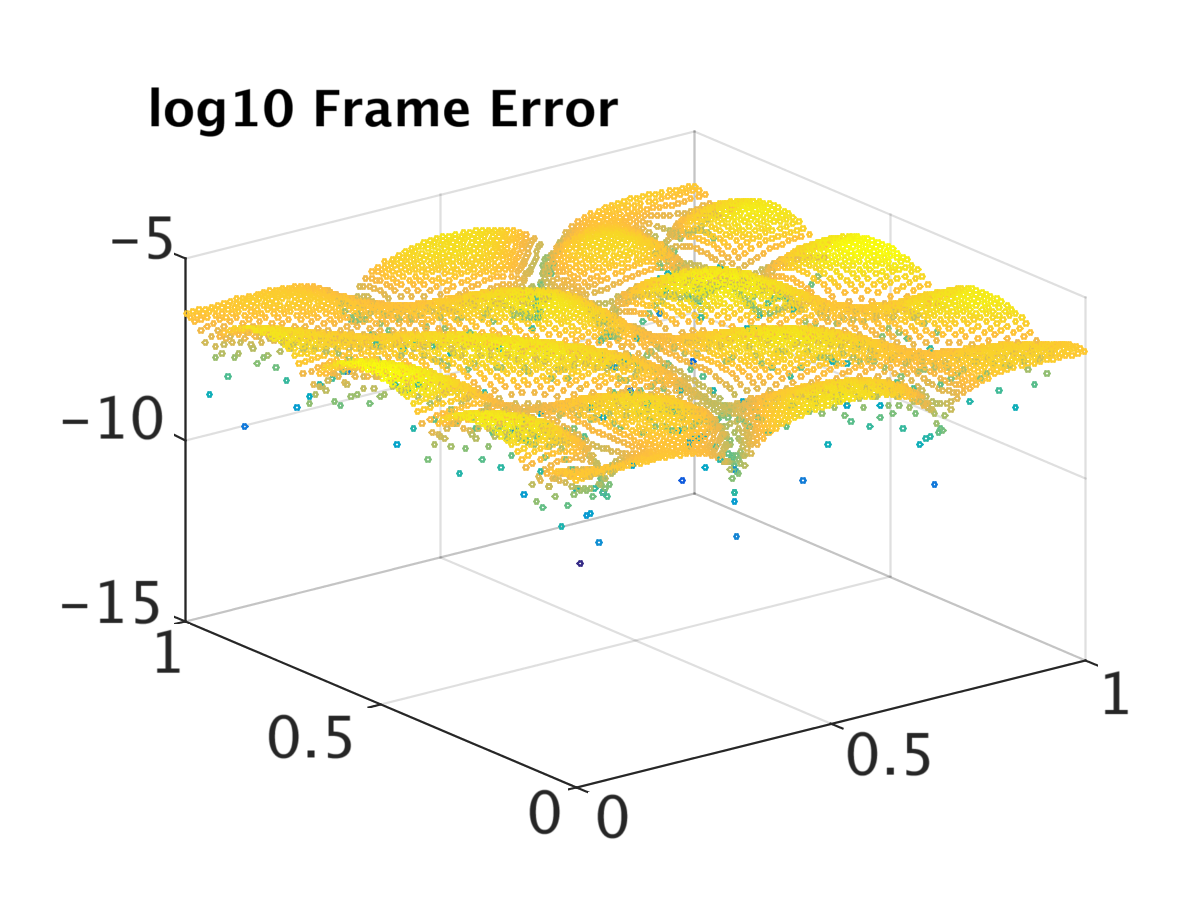}
\includegraphics[width=.245\textwidth]{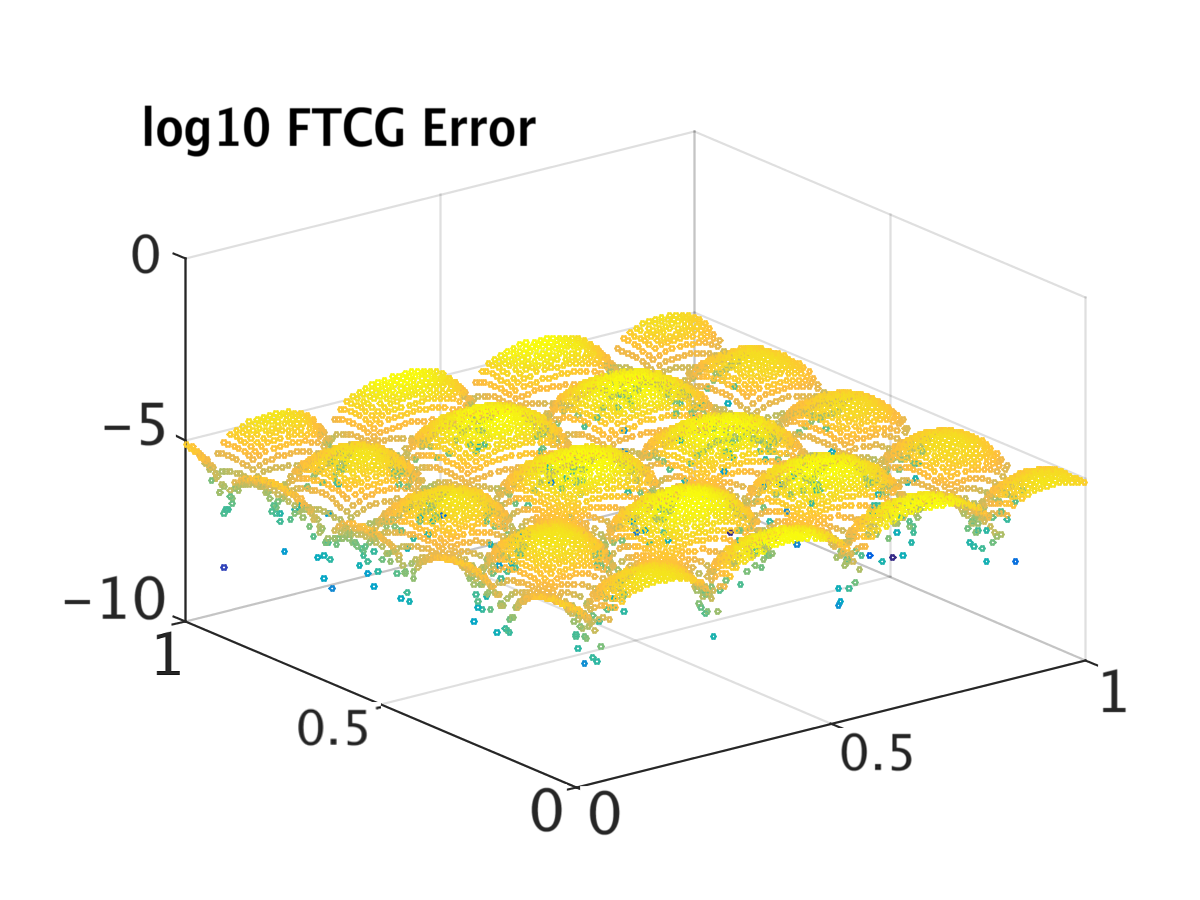}
\caption{Function reconstruction from asterisk data arrangement.}\label{fig:ast}
\end{minipage}

\begin{minipage}{1\textwidth}
\includegraphics[width=.5\textwidth]{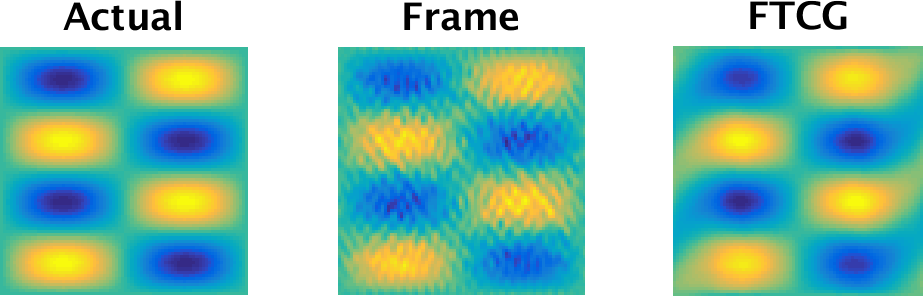}
\includegraphics[width=.245\textwidth]{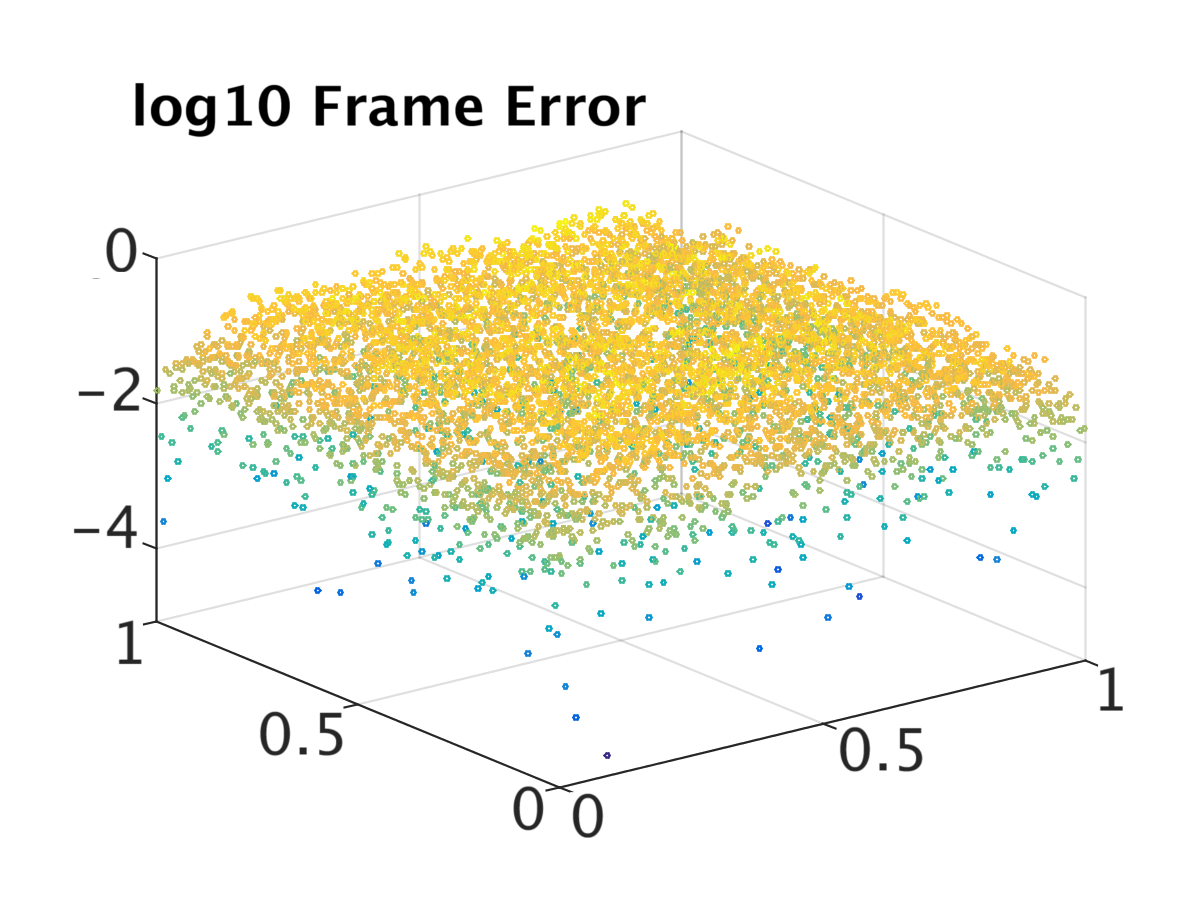}
\includegraphics[width=.245\textwidth]{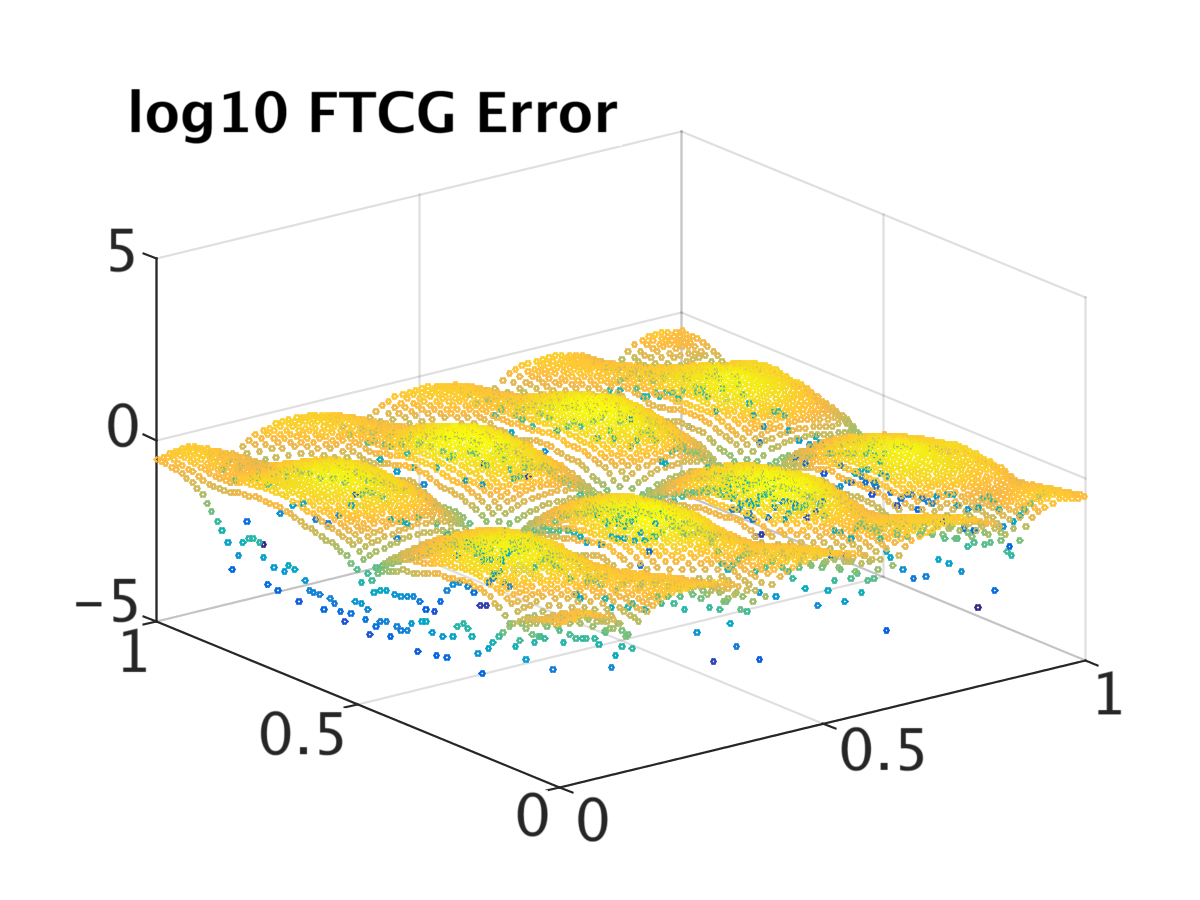}
\caption{Function reconstruction from SAS data arrangement.}\label{fig:sas}
\end{minipage}
\end{figure*}

%

In $1D$, it was proved in \cite{gelb2014frame} that $\psi_{n}$ is an admissible frame for $\phi_m$ if the raster is such that
\beqna 
\lambda_n\in [n-1/4,n+1/4].
\eeqna
Building on a $2D$ configuration that extends this optimal setting, as a first experiment we consider
\beqna
\bol\lambda_{\bol n}\in [n_1- 1/4,n_1+1/4]\times [n_2- 1/ 4, n_2+1/4]\,\,\forall\,\,\bol n\in\mathcal{I}_n.
\eeqna
In this case we used a band of 15 ($r=8$) for a matrix $T$ 
{which we define as
\begin{equation}
\label{tmatrix}
T= \Omega_2\Psi_2,
\end{equation}}
of dimension $900\times 900$ meaning that we used approximately 1.7\% of the information provided (and, therefore, 1.7\% of the computations).

Figure \ref{fig:nois} compares the frame and FTCG approximations for the example in (\ref{eq:eq2}).  In terms of PSNR, the frame approximation was able to achieve a value of 39.79dB while FTCG gave a value of 25.91dB.  It is important to note that no post processing was employed to reduce the effects of Fourier approximation on the the underlying function.\footnote{{More precisely, it is the the Fourier partial sum divided by the window function, $w(x)$.}}  Specifically, while the periodic function in (\ref{eq:eq2}) requires a significant number of Fourier coefficients to achieve a high resolution recovery, \cite{hesthaven2007spectral}, for piecewise smooth functions, the error approximation should be measured against the Fourier reconstruction (divided by the window function $w(x)$) of the underlying function. 

This experiment serves as a baseline proof of concept that the FTCG can produce satisfactory reconstructions at substantially less cost.  To test out more complex patterns that strain many of the frame theoretic assumptions, we present the next two experiments:  the asterisk and SAS-like arrangements.

\subsection{Asterisk Data Arrangement}

While one could think of the asterisk pattern given in Figure \ref{fig:jitgrid} as useful with regards to other tomographic imaging settings such as magnetic resonance imaging, a similar sampling could be found in circular SAS settings wherein a vehicle circles around a target, avoiding the shadowing present in side-scan \cite{callow2009circular}.  For this reason, and the fact that it is vastly different than the noisy grid, we sought out how the frame approximation and our 2D FTCG algorithm would perform.  We used $C_2$ with a band of 23 ($r=12$), representing 10\% of $T$ given in (\ref{tmatrix}).

\begin{figure}\centering
\includegraphics[width=.43\columnwidth]{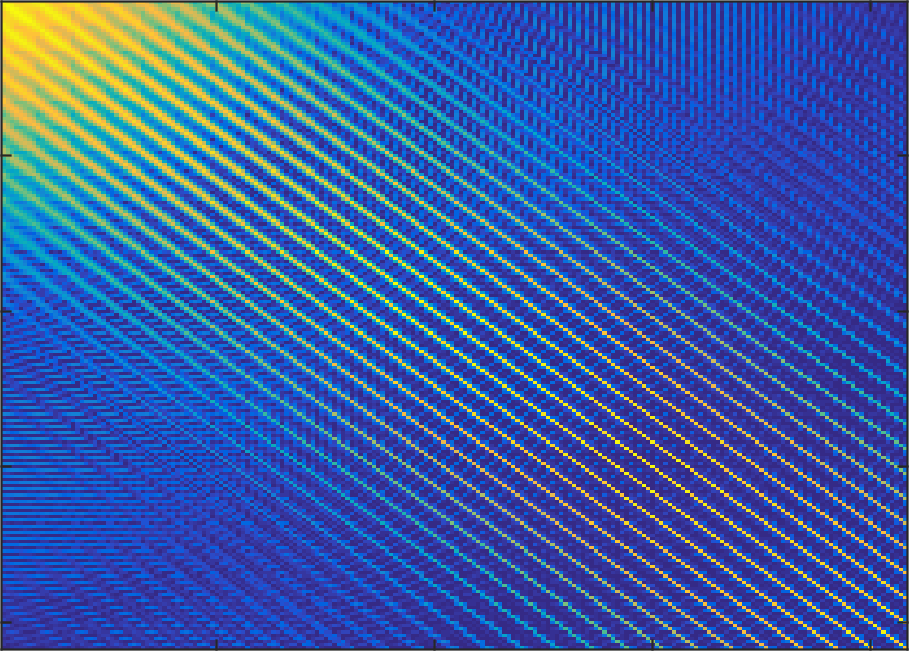}
\includegraphics[width=.43\columnwidth]{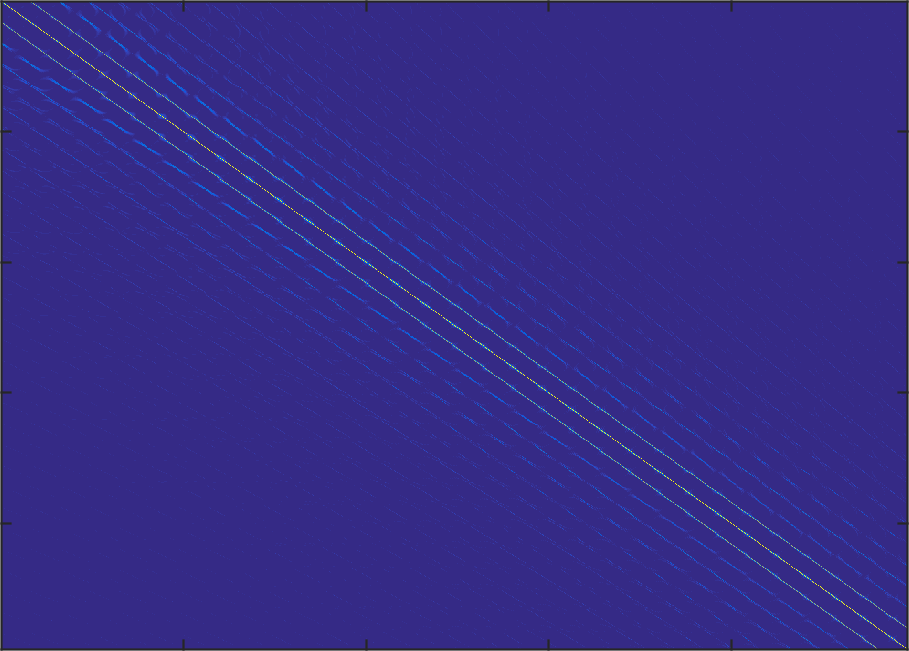}
\includegraphics[width=.056\columnwidth]{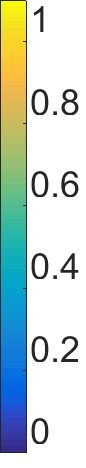}
\caption{Intensity maps of the magnitude values $T$ matrices {given by (\ref{tmatrix})} used for the asterisk sampling pattern (left) and side-scan SAS arrangement (right) experiments.}\label{fig:Tast}
\end{figure}

Both methods performed \emph{better} in our asterisk pattern case than the noisy grid (of which we had sound analytical inferences towards accuracy).  {This is likely due to the increased density in the low frequency domain, since our test function can be suitably resolved.}  In particular, observe the formation of matrix $T$ {given in (\ref{tmatrix})} in Figure \ref{fig:Tast} (left).   There is a high intensity region along the band (especially in the top left region) and decay off the center.  Given how well our images appear, the information contained on the off banded areas has minimal information and is, therefore, not really needed.  Figure \ref{fig:ast} shows the reconstructions and their $\log_{10}$ error.  The frame approximation yielded a PSNR value of 40.06dB and FTCG produced a value of 37.48dB meaning that, despite the banded structure and limited information, FTCG performed only marginally worse.  


\subsection{SAS Data Arrangement}

In our last experiment, we show a case where FTCG \emph{out performs} the frame approximation.  Here we use a side-scan-like arrangement (similar to the polar format common in synthetic aperture radar as well \cite{fan2014polar}).  The band was 15 which represented only $1.5\%$ of $T$.

While both reconstructions were challenged, the FTCG method was able to provide a clearer picture of the original function.  This fact is reflected by their PSNR values of 21.04 and 23.27 for the frame and FTCG approximations, respectively.  The point-wise error did lend some merit to the frame approximation over FTCG though this may reveal some issue with that metric.

A key reason why the FTCG yielded a higher PSNR may be shown in Figure \ref{fig:Tast} (right).  Observe the thin high intensity bands with large swaths of low magnitude information.  The diagonal itself is of the highest intensity and our method captures this while disregarding the rest.  {In particular we note that the condition number  $\kappa(\Psi) = 5851522$ while $\kappa(C_2) =  2.95$.} 
{Therefore we conclude that the FTCG serves as a pre-conditioner for the standard frame approximation.  We will further pursue using the FTCG as a pre-conditioner in future investigations.}  

\section{Discussion}
Understanding how to overcome the issues of $2D$ non-uniform Fourier data is essential to improving SAS.  While this approach to Sonar has dramatically improved underwater imaging, we see that most high-quality SAS reconstructions come from back-projection methods that are infeasible for computers without GPUs or other {very fast} processors.  To achieve higher resolutions in real time by automated underwater vehicles, methods have to be developed that {yield better and more consistent} performance than the standard $2D$ interpolations that typically cannot handle irregular motion and noise {due to ill conditioning}.  Hence we have presented the frame approximation as well as its less expensive surrogate, the  FTCG algorithm, as alternatives.  Both methods demonstrated improved resolution properties on synthetic data arranged on several rasters without any preprocessing to correct for imperfections.  Further work into adapting FTCG towards SAS should yield even better results with actual SAS data.

\bibliographystyle{IEEEbib}
\bibliography{refOceans2015}

\end{document}